# TWO EXAPLES OF EIGENFORS OF THE HODGE-DE RHAM LAPLACIAN

*Stepanov Sergey*

Department of Mathematics, Finance University under the Government of Russian Federation,
Leningradsky Prospect, 49-55, 125468 Moscow, Russian Federation
e-mail: *s.e.stepanov@mail.ru*

*Mikeš Josef*

Department of Algebra and Geometry, Palacky University,
77146 Olomouc, Listopadu 12, Czech Republic
e-mail: *josef.mikes@upol.cz*

*Tsyganok Irina*

Department of Mathematics, Finance University under the Government of Russian Federation,
Leningradsky Prospect, 49-55, 125468 Moscow, Russian Federation
e-mail: *i.i.tsyganok@mail.ru*

*Abstract*. In the present paper, we consider the Hodge-de Rham Laplacian $\Delta$ that acts on conformal Killing and projective Killing one-forms of a compact Riemannian manifold.

## 1. Introduction

Spectral geometry is a branch of mathematics that examines the relationship between the geometry of a manifold and the spectra of certain differential operators defined on this manifold. In the present paper, we consider the Hodge-de Rham Laplacian $\Delta$ that acts on conformal Killing and projective Killing one-forms of a compact Riemannian manifold. These forms are dual to well known conformal Killing and projective Killing vector fields, respectively.

It is known that the Lie algebra of group infinitesimal conformal transformations of an $n$-dimensional ($n \geq 3$) connected Riemannian manifold can be realized as the space of conformal Killing vector fields and has dimension $r \leq \frac{1}{2}(n+1)(n+2)$. On the other hand, group infinitesimal projective transformations of an $n$-dimensional ($n \geq 3$) connected Riemannian manifold can be realized as the space of projective Killing vector fields and has dimension $r \leq n(n+2)$. Hence, the corresponding eigenvalues of the Hodge-de Rham Laplacian $\Delta$ constitute two finite discrete spectrums. In the paper, we find lower and upper bounds of these spectrums.

## 2. The projective Killing form as an example of eigenforms of the Hodge-de Rham Laplacian

Let (*M, g*) be a connected Riemannian manifold of dimension $n \geq 2$ with the Levi-Civita connection $\nabla$. If a transformation of (*M, g*) preserves geodesics then it is called a *projective transformation*. Let $\xi$ be a vector field on (*M, g*) and $\{\varphi_t\}$ be a one-parameter group of local transformations generated by $\xi$. The vector field $\xi$ is called a *projective Killing vector field* or an *infinitesimal projective transformation* if its one-parameter group of $\{\varphi_t\}$ consists of local projective transformations (see [1, pp. 227-230]).

We define a *projective Killing form* $\omega$ by the condition $\omega(X) = g(\xi, X)$ for an infinitesimal projective transformation $\xi$ and an arbitrary smooth vector field $X$ on (*M, g*). This means that $\omega$ is dual to the vector field $\xi$ with respect to the metric *g*. The projective Killing form $\omega$ satisfies the following system of differential equations (see [1, p. 228])

$$(\nabla\Phi)(X,Y,Z) = 2\varphi(Z)g(X,Y) + \varphi(X)g(Z,Y) + \varphi(Y)g(X,Z) \qquad (2.1)$$

for $\Phi(X,Y) = (\nabla\omega)(X,Y) + (\nabla\omega)(Y,X)$ and $\varphi = -(n+1)^{-1} d\, d^*\omega$ where *d* is the operator of exterior differentiation and $d^*$ is the codifferentiation operator, which is the formal adjoint to *d*. We have the following theorem.

**Theorem 1.** *Let* (*M, g*) *be an n-dimensional* ($n \geq 2$) *compact and simply connected Riemannian manifold* (*M, g*) *which admits an infinitesimal projective transformation* $\xi$ *such that* $\Delta\omega = \lambda\omega$ *for the one-form* $\omega$ *dual to* $\xi$ *and a real positive number* $\lambda$. *Suppose that* $\rho$ *and P are the smallest and the largest eigenvalues of the matrix* $\|Ric\|$ *on* (*M, g*) *then* $2\rho \leq \lambda \leq 2(n+1)^{-1}(n-1)P$. *If the equality* $\lambda = 2\rho$ *is attained then* $\xi$ *is a Killing vector field. On the other hand, if the equality* $\lambda = 2(n+1)^{-1}(n-1)P$ *is attained then* $\xi = \text{grad } f$ *for some smooth scalar function* $f : M \to \mathbb{R}$. *Moreover, if the condition* $\Delta f = 2(n+1)\alpha f$ *is satisfied, then* (*M, g*) *is isometric to a Euclidian sphere of radius* $1/\sqrt{\alpha}$.

**Proof.** Yano has proved (see [2, p. 45]) that an arbitrary projective Killing one-form $\omega$ satisfies the equation

$$\Delta\omega = 2Ric^*(\xi) + 2(n+1)^{-1} d\, d^*\omega \qquad (2.2)$$

where $\Delta = d\, d^* + d^*d$ is the well known Hodge-de Rham Laplacian and $Ric^*$ is a symmetric algebraic operator which is defined by the following condition $Ric(\xi, X) = g(Ric^*(\xi), X)$ for the Ricci tensor *Ric* and an arbitrary vector field $X \in C^\infty TM$. It is obvious that the equation (2.2) can be rewritten in the form

$$\Delta\omega = 2(n+1)(n-1)^{-1} Ric^*(\xi) + 2(n-1)^{-1} d^*d\, \omega. \qquad (2.3)$$

Let a positive number $\lambda$ be an eigenvalue of the Hodge-de Rham Laplacian $\Delta$ such that the corresponding eigenform one-form $\omega$ is a projective Killing form, i.e. $\Delta\omega = \lambda\omega$. If we assume that (M, g) is compact and the Ricci tensor Ric is positive then there exist two positive numbers $\rho$ and $P$ which are the smallest and the largest eigenvalues of the matrix $\|Ric\|$ on (M, g). In this case we have the following inequalities $\rho g(X,X) \leq Ric(X,X) \leq P g(X,X)$ for an arbitrary vector field $X \in C^\infty TM$. Further, from these inequalities and the equations (2.2) and (2.3) we can obtain the following inequalities $2\rho \leq \lambda \leq 2(n+1)^{-1}(n-1)P$.

If we suppose that the equality $\lambda = 2\rho$ is attained then from (2.1) we obtain that $d^*\omega = 0$, and hence $\xi$ is a *Killing vector field* (see [2, p. 46]). On the other hand, if the equality $\lambda = (n+1)^{-1}(n-1)P$ is attained then from (2.3) we obtain that $d\omega = 0$. We recall here that the first Betti number $b_1(M) = 0$ for a compact manifold (M, g) with the positive definite Ricci tensor Ric (see [2, p. 42]). In this case, Ker $d$ = Im $d$, and hence $\omega = df$ for some smooth scalar function $f: M \to \mathbb{R}$. At the same time from the equations (2.1), we conclude that the function $f$ satisfies the following system of differential equations of order three

$$2(\nabla^3 f)(X,Y,Z) = 2\varphi(Z)g(X,Y) + \varphi(X)g(Z,Y) + \varphi(Y)g(X,Z) \tag{2.4}$$

where $\varphi = -(n+1)^{-1} d\Delta f$. For the case $\Delta f = \mu f$ where $\mu$ is a positive constant we can rewrite (2.4) in the form

$$(\nabla^3 f)(X,Y,Z) + \alpha\left(2df(Z)g(X,Y) + df(X)g(Z,Y) + df(Y)g(X,Z)\right) = 0. \tag{2.5}$$

where $\alpha = \frac{1}{2}\mu(n+1)^{-1}$. Simultaneously, we recall (see [3]) that a complete connected and simply connected Riemannian manifold (M, g) of dimension $n \geq 2$ which admits a non-constant solution $f$ satisfying the differential equation (2.5) is globally isometric to a sphere ($\mathbb{S}^n$, $g_0$) of radius $1/\sqrt{\alpha}$ in the Euclidian space $\mathbb{R}^{n+1}$. The theorem is proved.

**Remark**. It is well known that the second eigenvalue $\mu_2$ of the Laplacian $\Delta = d^*\nabla$ which acts on smooth scalar functions $f: \mathbb{S}^n \to \mathbb{R}$ for an arbitrary sphere ($\mathbb{S}^n$, $g_0$) of radius $1/\sqrt{\alpha}$ is $\mu_2 = 2(n+1)\alpha$ and each eigenfunction $f$ corresponding to $\mu_2 = 2(n+1)\alpha$ satisfies (2.5). In our case from the system of differential equations (2.5) we also deduce that $\mu = 2(n+1)\alpha$.

### 3. The conformal Killing form as an example of eigenforms of the Hodge-de Rham Laplacian

Let (*M, g*) be a connected Riemannian manifold of dimension $n \geq 3$. *Conformal Killing vector fields* can be considered as a natural generalization of Killing vector fields. They are also called *infinitesimal conformal transformations* because any conformal Killing vector $\xi$ generates a local one-parameter group $\{\varphi_t\}$ of local conformal transformations of (*M, g*).

We define a *conformal Killing form* $\omega$ by the condition $\omega(X) = g(\xi, X)$ for an infinitesimal conformal transformation $\xi$ and an arbitrary smooth vector field $X$ on (*M, g*). The conformal Killing form $\omega$ satisfies the following system of differential equations (see [1, p. 228])

$$(\nabla \omega)(X,Y) + (\nabla \omega)(Y,X) + 2n^{-1} d^* \omega\, g(X,Y) = 0 \tag{3.1}$$

for arbitrary smooth vector fields *X* and *Y* on (*M, g*). We have the following theorem.

**Theorem 2.** *Let* (*M, g*) *be an n-dimensional* ($n \geq 2$) *compact and simply connected Riemannian manifold* (*M, g*) *which admits an infinitesimal conformal transformation* $\xi$ *such that* $\Delta \omega = \lambda \omega$ *for the one-form* $\omega$ *dual to* $\xi$ *and a real positive number* $\lambda$. *Suppose that* $\rho$ *and P are the smallest and the largest eigenvalues of the matrix* $\|Ric\|$ *on* (*M, g*) *then* $n(n-1)^{-1}\rho \leq \lambda \leq 2P$. *If the equality* $\lambda = n(n-1)^{-1}\rho$ *is attained then* $\xi$ *is a Killing vector field. On the other hand, if the equality* $\lambda = 2P$ *is attained, then* (*M, g*) *is conformally diffeomorphic an Euclidian sphere* ($\mathbb{S}^n$, $g_0$). *Moreover, if the condition* $\Delta f = n\alpha f$ *is satisfied, then* (*M, g*) *is isometric to an Euclidian sphere of radius* $1/\sqrt{\alpha}$.

**Proof**. Lichnerowicz has shown (see [2, p. 47]) that a necessary and sufficient condition for $\xi$ to be a *conformal Killing vector field* on a compact Riemannian manifold (*M, g*) is

$$\Delta \omega = 2 Ric^*(\xi) - (1 - 2/n)\, d\, d^* \omega = 0 \tag{3.2}$$

for the 1-form $\omega$ dual to the vector field $\xi$ with respect to the metric *g*. The 1-form $\omega$ will be called *conformal Killing form*. It is obvious that the equation (3.2) can be rewritten in the form

$$\Delta \omega = n(n-1)^{-1} Ric^*(\xi) + \tfrac{1}{2}(n-1)^{-1}(n-2) d^* d\, \omega. \tag{3.3}$$

Let a positive number $\lambda$ be an eigenvalue of the Hodge-de Rham Laplacian $\Delta$ such that the corresponding eigenform one-form $\omega$ is a conformal Killing form, i.e. $\Delta \omega = \lambda \omega$. We assume that (*M, g*) is compact and the Ricci tensor *Ric* is positive. In this case there exist two positive numbers $\rho$ and *P* which are the smallest and the largest eigenvalues of the matrix $\|Ric\|$ on (*M, g*). In this case we have two inequalities $\rho g(X,X) \leq Ric(X,X) \leq P g(X,X)$ for an arbitrary vector field $X \in C^\infty TM$. Further, from these inequalities and the equations (3.2) and (3.3) we can obtain the inequalities $n(n-1)^{-1}\rho \leq \lambda \leq 2P$.

If we suppose that the equality $\lambda = n(n-1)^{-1}\rho$ is attained then from (3.2) we obtain that $d^*\omega = 0$, and hence $\xi$ is a *Killing vector field* (see [2, p. 46]). On the other hand, if the equality $\lambda = 2P$ is attained then from (2.3) we obtain that $d\omega = 0$. At the same time from the equations (3.1), we conclude that the function $f$ satisfies the following system of differential equations of order two

$$(\nabla^2 f)(X,Y) + n^{-1}\Delta f\, g(X,Y) = 0 \tag{3.4}$$

for arbitrary smooth vector fields $X$ and $Y$ on $(M, g)$. According to [4], the existence of such a non-constant smooth function $f$ means that the compact manifold $(M, g)$ is conformally diffeomorphic to a sphere $(\mathbb{S}^n, g_0)$ in the Euclidian space $\mathbb{R}^{n+1}$.

For the case $\Delta f = \mu f$ where $\mu$ is a positive constant we can rewrite (3.4) in the form

$$(\nabla^2 f)(X,Y) + n^{-1}\mu f \cdot g(X,Y) = 0. \tag{3.5}$$

In this case $(M, g)$ is globally isometric to a sphere $(\mathbb{S}^n, g_0)$ of radius $\sqrt{n/\mu}$ in the Euclidian space $\mathbb{R}^{n+1}$ (see [5]).

**Remark**. We recall that the first eigenvalue $\mu_1$ of the Laplacian $\Delta = d^*\nabla$ which acts on smooth scalar functions $f : \mathbb{S}^n \to \mathbb{R}$ for an arbitrary sphere $(\mathbb{S}^n, g_0)$ of radius $1/\sqrt{\alpha}$ is $\mu_1 = n\alpha$ and each eigenfunction $f$ corresponding to $\mu_1 = n\alpha$ satisfies the following system of differential equations $\nabla^2 f + \alpha f \cdot g_0 = 0$. Simultaneously, from (3.5) we deduce that $\mu = n\alpha$.

## 4. Appendix

Tanno's rigidity theorem (which has various important geometric applications) as state below is well-known (see [3]).

**Theorem 3**. *Let $(M, g)$ be a complete and simple and simply connected Riemannian manifold. If $(M, g)$ admits a non-constant function f satisfying the following system of differential equations of order three*

$$(\nabla^3 f)(X,Y,Z) + k\,(2df(Z)g(X,Y) + df(X)g(Z,Y) + df(Y)g(X,Z)) = 0. \tag{4.1}$$

*for some positive constant k and arbitrary smooth vector fields X, Y, Z, it is necessary and sufficient that $(M, g)$ is isometric to a Euclidian sphere $\mathbb{S}^n$ of radius $1/\sqrt{k}$.*

Now we transform (4.1) to the following system of differential equations

$$(\nabla^3 f)(X,Y,Z) + 2\varphi(Z)g(X,Y) + \varphi(Z)(X)g(Z,Y) + \varphi(Z)(Y)g(X,Z) = 0 \tag{4.2}$$

for $\varphi = \frac{1}{2(n+1)} d \Delta f$. Next, we express the assumption that the following generalized Tanno theorem holds.

**Theorem 4**. *Let (M, g) be a complete and simple and simply connected Riemannian manifold. If (M, g) admits a non-constant function f satisfying* (4.2), *then (M, g) is conformally diffeomorphic to a Euclidian sphere* ($\mathbb{S}^n$, $g_0$).

The generalized Tanno theorem is consistent with the statement of the Tashiro theorem (see [4]). We note that the proof of this theorem we do not know, but in this case our Theorem 1 can be rewritten in the following complete form.

**Theorem 5.** *Let (M, g) be an n-dimensional ($n \geq 2$) compact and simply connected Riemannian manifold (M, g) which admits an infinitesimal projective transformation $\xi$ such that $\Delta \omega = \lambda \omega$ for the one-form $\omega$ dual to $\xi$ and a real positive number $\lambda$. Suppose that $\rho$ and P are the smallest and the largest eigenvalues of the matrix $\|Ric\|$ on (M, g) then $2\rho \leq \lambda \leq 2(n+1)^{-1}(n-1)P$. If the equality $\lambda = 2\rho$ is attained then $\xi$ is a Killing vector field. On the other hand, if the equality $\lambda = 2(n+1)^{-1}(n-1)P$ is attained then (M, g) is conformally diffeomorphic to a Euclidian sphere* ($\mathbb{S}^n$, $g_0$). *Moreover, if the condition $\Delta f = 2(n+1)\alpha f$ is satisfied, then (M, g) is isometric to a Euclidian sphere of radius $1/\sqrt{\alpha}$.*

The new formulation of our Theorem 1 is consistent with the statement of our Theorem 2.